\documentclass[twoside,a4paper,12pt,centertags,reqno]{amsart} 
\usepackage{amsmath,amssymb,verbatim,vmargin}
\usepackage{color}
\usepackage{tikz}
\usepackage{hyperref}
\hypersetup{colorlinks,bookmarks=true,linktocpage=true,citecolor=blue,linkcolor=magenta}

\usepackage{eulervm}   

\allowdisplaybreaks 
\usepackage{color}
\usepackage[all]{xy} 

\usepackage{mathtools}
\usepackage{MnSymbol} 

\theoremstyle{plain}
\newtheorem{thm}{Theorem}[section]

\theoremstyle{definition}

\newtheorem{rem}[thm]{Remark}

\usepackage{enumerate}

\newcommand{\bN}{{\mathbb N}}

\def\barint_#1{\mathchoice
            {\mathop{\vrule width 6pt
height 3 pt depth -2.5pt
                    \kern -9.5pt
\intop \kern -4pt}\nolimits_{#1}}%
            {\mathop{\vrule width 5pt height
3 pt depth -2.6pt
                    \kern -6.5pt
\intop \kern -4pt}\nolimits_{#1}}%
            {\mathop{\vrule width 5pt height
3 pt depth -2.6pt
                    \kern -6pt
\intop \kern -4pt}\nolimits_{#1}}%
            {\mathop{\vrule width 5pt height
3 pt depth -2.6pt
          \kern -6pt \intop \kern -4pt}\nolimits_{#1}}}
          
           \def\bariint_#1{\mathchoice
            {\mathop{\vrule width 15pt
height 3 pt depth -2.5pt
                    \kern -15.8pt
\intop \kern -8pt\intop \kern -4pt}\nolimits_{#1}}%
            {\mathop{\vrule width 9pt height
3 pt depth -2.6pt
                    \kern -10.5pt
\intop \kern -8pt\intop \kern -4pt}\nolimits_{#1}}%
            {\mathop{\vrule width 9pt height
3 pt depth -2.6pt
                    \kern -10pt
\intop \kern -8pt\intop \kern -4pt}\nolimits_{#1}}%
            {\mathop{\vrule width 9pt height
3 pt depth -2.6pt
          \kern -8pt \intop \kern -10pt\intop \kern -4pt}
      \nolimits_{  #1}}}

\def\barintlim_#1{\mathchoice
            {\mathop{\vrule width 6pt
height 3 pt depth -2.5pt
                    \kern -8.8pt
\intop \kern -4pt}\limits_{#1}}%
            {\mathop{\vrule width 5pt height
3 pt depth -2.6pt
                    \kern -6.5pt
\intop \kern -4pt}\limits_{#1}}%
            {\mathop{\vrule width 5pt height
3 pt depth -2.6pt
                    \kern -6pt
\intop \kern -4pt}\limits_{#1}}%
            {\mathop{\vrule width 5pt height
3 pt depth -2.6pt
          \kern -6pt \intop \kern -4pt}\limits_{#1}}}
          
           \def\bariintlim_#1{\mathchoice
            {\mathop{\vrule width 15pt
height 3 pt depth -2.5pt
                    \kern -15.8pt
\intop \kern -8pt\intop \kern -4pt}\limits_{#1}}%
            {\mathop{\vrule width 9pt height
3 pt depth -2.6pt
                    \kern -10.5pt
\intop \kern -8pt\intop \kern -4pt}\limits_{#1}}%
            {\mathop{\vrule width 9pt height
3 pt depth -2.6pt
                    \kern -10pt
\intop \kern -8pt\intop \kern -4pt}\limits_{#1}}%
            {\mathop{\vrule width 9pt height
3 pt depth -2.6pt
          \kern -8pt \intop \kern -10pt\intop \kern -4pt}
      \limits_{  #1}}}
          
\renewcommand{\iint}{\int \kern -3pt\int}       

\numberwithin{equation}{section}
\usepackage[symbol]{footmisc}

\title{A proof of the Prelov conjecture}

\author{Yi C. Huang} 
\address{School of Mathematical Sciences, Nanjing Normal University, Nanjing 210023, People's Republic of China}
\email{Yi.Huang.Analysis@gmail.com}
\urladdr{https://orcid.org/0000-0002-1297-7674}

\author{Fei Xue} 
\address{School of Mathematical Sciences, Nanjing Normal University, Nanjing 210023, People's Republic of China}
\email{05429@njnu.edu.cn}

\date{\today} 

\thanks{Research of the authors is partially supported by the National NSF grants of China (nos. 11801274 and 12201307) and the Jiangsu Provincial NSF grant (no. BK20210555).
This paper is started while YCH is on leave, funded by CSC Postdoctoral/Visiting Scholar Program (no. 202006865011), at LAGA of Universit\'e Sorbonne Paris Nord.
YCH would also like to thank Profs. X.-Y. Gui and Y. Yang for initial helpful communications in classical information theory.}

\begin{document}

\begin{abstract}
In this paper we present a complete proof of a conjecture due to V. V. Prelov in 2010 about an information inequality for the binary entropy function.
\end{abstract}

\maketitle


\section{Introduction}

Consider the binary entropy function 
\begin{equation} \label{eqn:h}
h(x)=-x\log x-(1-x)\log(1-x),\quad x\in[0,1].
\end{equation}
In computing the maximum of mutual information of several random variables via the variational distance 
(see Pinsker \cite{Pin05} and Prelov \cite{Pre09}),
Prelov established in \cite{Pre10} some new inequalities for the binary entropy function.

\begin{thm}[Prelov, 2010]
For $p>q$, we have
\begin{equation} \label{eqn:Prelov}
p^n(1+q^n-p^n)h(q^n)>q^n(1+p^n-q^n)h(p^n),
\end{equation}
if 
$$0<q<q_0,\quad \text{or}\,\, \,\,q_1<q<1/2,\quad \text{or}\,\, \,\,n>n_0,$$ 
where $q_0$, $q_1$ and $n_0$ are some positive numbers.
\end{thm}

For $n\in\{1,\infty\}$, or for $(p,q)\in\{(1,0), (1/2,1/2)\}$, \eqref{eqn:Prelov} is an equality.
Prelov conjectures that \eqref{eqn:Prelov} is valid for all $q\in(0,1/2)$, $p=1-q$ and all integers $n\geq2$.

\begin{thm} \label{thm:Prelov2}
The above conjecture of Prelov is true.
\end{thm}

For the impact of this result to information theory, see \cite[Proposition 1]{Pre10}.

\section{Proof of Theorem \ref{thm:Prelov2} in quadratic case}

Using $p=1-q$ we see that \eqref{eqn:Prelov} for $n=2$ is equivalent to
\begin{equation} \label{eqn:Prelov'}
ph(q^2)>qh(p^2).
\end{equation}
Using the definition \eqref{eqn:h} and $p+q=1$, we have further equivalences
\begin{equation} \label{eqn:Prelov''}
\begin{aligned}
\eqref{eqn:Prelov'}&\Longleftrightarrow -2pq^2\log q-p^2(1+q)\log(1+q)-p^2(1+q)\log p\\
&\qquad\qquad>-2qp^2\log p-q^2(1+p)\log(1+p)-q^2(1+p)\log q\\
&\Longleftrightarrow q^2\left[(1-p)\log(1-p)+(1+p)\log(1+p)\right]\\
&\qquad\qquad>p^2\left[(1-q)\log(1-q)+(1+q)\log(1+q)\right].\\
\end{aligned}
\end{equation}
To this end, we introduce
$$\lambda(x)=\frac{\left[(1-x)\log(1-x)+(1+x)\log(1+x)\right]}{x^2}$$ 
and
$$\lambda'(x)=\frac{(x-2)\log(1-x)-(x+2)\log(1+x)}{x^3}=:\frac{\beta(x)}{x^3}.$$
Note that
$$\gamma(x):=\beta'(x)=\log(1-x)+\frac{x-2}{x-1}-\log(1+x)-\frac{x+2}{x+1}$$
and
$$\gamma'(x)=\frac{x}{(x-1)^2}-\frac{x}{(x+1)^2}=\frac{4x^2}{(1-x^2)^2}>0, \quad x\in(0,1).$$
Now, since $\gamma(0)=0$, so $\gamma>0$ on $(0,1)$ and $\beta$ is strictly increasing.
Similarly, since $\beta(0)=0$, so $\beta>0$ on $(0,1)$ and $\lambda$ is strictly increasing.
In particular, 
$$p>\frac12>q\Longrightarrow\lambda(p)>\lambda(q).$$ 
This together with the equivalence \eqref{eqn:Prelov''} proves the theorem.

\begin{rem}
Observe that we have the simple inequality
\begin{equation} \label{eqn:elem}
h(p^2)<2ph(p),\quad 0<p<1.
\end{equation}
Indeed, by the definition of $h$, \eqref{eqn:elem} is equivalent to
$$g(p)=(1+p)\log(1+p)+(1-p)\log(1-p)>0,\quad 0<p<1,$$
which is obvious by $g(0)=0$ and $g'(p)>0$ for $0<p<1$.
Now, using \eqref{eqn:Prelov'} and \eqref{eqn:elem},
$$h(q^2)+h(p^2)<\frac{h(q^2)}{q}<2h(q)=2h(p),\quad 0<p=1-q<1.$$
This gives an alternative (non-information-theoretic) proof of \cite[(10)]{Pre10} for $n=2$.
\end{rem}

\section{Proof of Theorem \ref{thm:Prelov2} in general case}

Given Theorem \ref{thm:Prelov2} in quadratic case (already verified in previous section), 
the idea is then to tailor for \eqref{eqn:Prelov} the following symmetric function
$$F(x,y)=\frac{1}{x-y}\left(\frac{1+y-x}{y}h(y)-\frac{1+x-y}{x}h(x)\right)$$ 
defined on the off-diagonal square $S':=\{x\neq y: 0<x,y<1\}$.
So we have verified
$$F(x,y)>0\quad\text{on}\quad \Gamma_2:=\left\{(x,y)\in S': x^{\frac12}+y^{\frac12}=1\right\},$$
and we want to show, for any $3\leq n\in\bN$, 
$$F(x,y)>0\quad\text{on}\quad  \Gamma_n:=\left\{(x,y)\in S': x^{\frac1n}+y^{\frac1n}=1\right\}.$$
Note that $\Gamma_n$ lies below $\Gamma_2$.  
Now, the trick is to show $F(x,y)$ is decreasing along the diagonal direction $(1,1)$.
By direct computations this directional derivative is
$$\begin{aligned}G(x,y)&:=\frac{d}{dt}F(x+t,y+t)\bigg|_{t=0}\\
&=\frac{1}{x-y}\left(\frac{1+y-x}{y}h'(y)-\frac{1+y-x}{y^2}h(y)-\frac{1+x-y}{x}h'(x)+\frac{1+x-y}{x^2}h(x)\right)\\
&=\frac{g(x)-g(y)}{x-y}+g(x)+g(y),
\end{aligned}$$
where
$$g(x)=\frac{h(x)-xh'(x)}{x^2}=\frac{2x-1}{x^2}\log(1-x).$$
Note that $g(1/2)=0$, and by routine computations, we have
\begin{equation} \label{eqn:g}
g'<0\quad \text{on}\quad (0,1); \quad g''>0\quad \text{on}\quad  (0,1/2).
\end{equation}
By symmetry we can consider $ G$ only for $x<y$. We then argue in two cases.

\textbf{Case I.} If $x<y$ and $y>1/2$, by monotonicity and $1+\frac{1}{x-y}<0$,
$$\begin{aligned}
 G(x,y)&=\left(1+\frac{1}{x-y}\right)g(x)+\left(1-\frac{1}{x-y}\right)g(y)\\
&<\left(1+\frac{1}{x-y}\right)g(y)+\left(1-\frac{1}{x-y}\right)g(y)=2g(y)<0.
\end{aligned}$$

\textbf{Case II.} If $x<y$ and $y\leq1/2$, by convexity and monotonicity,
$$ G(x,y)<\frac{g(x)-g(1/2)}{x-1/2}+g(x)+g(x)=\frac{4}{x}\ln(1-x)<0.$$

\noindent Theorem \ref{thm:Prelov2} in general case is proved.

\begin{rem}
For completeness, we also justify \eqref{eqn:g} with details. Indeed,
$$\begin{aligned}
g'(x)&=\frac{2\left(1-x\right)^{2}\ln\left(1-x\right)-2x^{2}+x}{x^{3}\left(1-x\right)}\\
&<\frac{-2(1-x)^2 x-2x^2+x}{x^3(1-x)}=-\frac{1-2x+2x^2}{x^2(1-x)}<0.
\end{aligned}$$
$$g''(x)=-\frac{x(6x^2-9x+4)+(6-4x)(1-x)^2\ln(1-x)}{x^4(1-x)^2}=:-\frac{a(x)}{x^4(1-x)^2}.$$
To see $a(x)<0$ when $x\in(0,1/2)$, note that $a(0)=0$, $a(1/2)<0$, and
$$a''(x)=4(2+4x+(7-6x)\ln(1-x))>4(2+4x-(7-6x)\times 1.4x)>0.$$
\end{rem}

\bigskip

\section*{\textbf{Compliance with ethical standards}}

\bigskip

\textbf{Conflict of interest} The authors have no known competing financial interests
or personal relationships that could have appeared to influence this reported work.

\bigskip

\textbf{Availability of data and material} Not applicable.

\bigskip

\bibliographystyle{alpha}

\bibliography{HuaY-XueF-PrelovComplete} 
 
\end{document}